\documentclass{emsprocart}
\contact[geisser@usc.edu]{Thomas Geisser, Dep. of Math. University of 
Southern California, 3620 S Vermont Av. KAP 108, Los Angeles, CA 90089, USA}

\usepackage{latexsym}
\usepackage{amsfonts}
\usepackage{amsmath}
\usepackage{amssymb}
\usepackage{amscd}
\usepackage{enumerate}
    \newtheorem{theorem}                    {Theorem}       [section]
    \newtheorem{lemma}      [theorem]       {Lemma}
    \newtheorem{corollary}  [theorem]       {Corollary}
    \newtheorem{proposition}[theorem]       {Proposition}

\catcode`@=11
\atdef@ I#1I#2I{\CD@check{I..I..I}{\llap{$\m@th\vcenter{\hbox
  {$\scriptstyle#1$}}$}
  \rlap{$\m@th\vcenter{\hbox{$\scriptstyle#2$}}$}&&}}
\atdef@ E#1E#2E{\ampersand@
  \ifCD@ \global\bigaw@\minCDarrowwidth \else \global\bigaw@\minaw@ \fi
  \setboxz@h{$\m@th\scriptstyle\;\;{#1}\;$}%
  \ifdim\wdz@>\bigaw@ \global\bigaw@\wdz@ \fi
  \@ifnotempty{#2}{\setbox@ne\hbox{$\m@th\scriptstyle\;\;{#2}\;$}%
    \ifdim\wd@ne>\bigaw@ \global\bigaw@\wd@ne \fi}%
  \ifCD@\enskip\fi
    \mathrel{\mathop{\hbox to\bigaw@{}}%
      \limits^{#1}\@ifnotempty{#2}{_{#2}}}%
  \ifCD@\enskip\fi \ampersand@}
\catcode`@=\active

\newcommand{\chr}{\operatorname{char}}

\newcommand{\Hom}{\operatorname{Hom}}
\newcommand{\End}{\operatorname{End}}

\newcommand{\Gal}{\operatorname{Gal}}

\newcommand{\trdeg}{\operatorname{trdeg}}

\renewcommand{\lim}{\operatornamewithlimits{lim}}
\newcommand{\colim}{\operatornamewithlimits{colim}}

\newcommand{\Z}{{{\mathbb Z}}}
\newcommand{\Q}{{{\mathbb Q}}}

\newcommand{\bm}{{c}}
\newcommand{\B}{{\mathbb Z}^c}

\newcommand{\F}{{{\mathbb F}}}

\newcommand{\et}{{\text{\rm et}}}

\newcommand{\C}{{\cal C}}

\newcommand{\proofend}{\hfill $\Box$ \\}
\newcommand{\rem}{\noindent {\it Remark. }}

\title[Parshin's conjecture revisited]{Parshin's conjecture revisited}
\author[Thomas Geisser]{Thomas Geisser\thanks{Supported in part by NSF grant No.0556263}}

\begin{document}


\maketitle

\section{Introduction}
Parshin's conjecture states that $K_i(X)_\Q=0$ for $i>0$ and $X$
smooth and projective over a finite field $\F_q$. The purpose of
this paper is to break up Parshin's conjecture into several
independent statements, in the hope that each of them is easier to
attack individually. If $CH_n(X,i)$ is Bloch's higher Chow group of
cycles of relative dimension $n$, then in view of $K_i(X)_\Q\cong
\oplus_n CH_n(X,i)_\Q$, Parshin's conjecture is equivalent to
Conjecture $P(n)$ for all $n$, stating that $CH_n(X,i)_\Q=0$ for
$i>0$, and all smooth and projective $X$. We show assuming resolution
of singularities that
Conjecture $P(n)$ is equivalent to the conjunction of three
conjectures $A(n)$, $B(n)$ and $C(n)$, and give several equivalent
versions of these conjectures. This is most conveniently formulated
in terms of weight homology. We define $H^W_*(X,\Q(n))$ to be the
homology of the complex $CH_n(W(X))_\Q$, where $W(X)$ is the
weight complex defined by Gillet-Soul\'e \cite{gilletsoule} shifted by $2n$.
Then, in a nutshell, Conjecture $A(n)$ states
that for all schemes $X$ over $\F_q$, the niveau spectral sequence
of $CH_n(X,*)_\Q$ degenerates to one line, Conjecture $C(n)$ states
that for all schemes $X$ over $\F_q$ the niveau spectral sequence of
$H_*^W(X,\Q(n))$ degenerates to one line, and Conjecture $B(n)$
states that, for all $X$ over $\F_q$, the two lines are isomorphic.
The conjunction of $A(n)$, $B(n)$, and $C(n)$ clearly implies
$P(n)$, because $H_i^W(X,\Q(n))=0$ for $i\not=2n$ and $X$ smooth and
projective over $\F_q$, and we show the converse.
Note that in the above formulation, Parshin's conjecture implies
statements on higher Chow groups not only for smooth and projective
schemes, but gives a way to calculate $CH_n(X,i)_\Q$ for all $X$.

A more concrete version of $A(n)$ is that, for every smooth and projective
scheme $X$ of dimension $d$ over $\F_q$, $CH_n(X,i)_\Q=0$ for $i>d-n$.
A reformulation of $B(n)$ is that for every smooth and projective
scheme $X$ of dimension $d>n+1$, the following sequence is exact
$$0\to K_{d-n}^M(k(X))_\Q\to\oplus_{x\in X^{(1)}}K_{d-n-1}^M(k(x))_\Q
\to \oplus_{x\in X^{(2)}}K_{d-n-2}^M(k(x))_\Q$$
and that this sequence is exact at $K_1^M(k(X))_\Q$ if $d=\dim X=n+1$.

In the second half of the paper, we focus on the case $n=0$,
because of its applications in \cite{ichkato}.
A different version of weight homology has been
studied by Jannsen \cite{jannsen}, and he proved that Conjecture $C(0)$
holds under resolution of singularities. We use this to
give two more versions of Conjecture $P(0)$. The first is that
there is an isomorphism from higher Chow groups with $\Q_l$-coefficients
to $l$-adic cohomology, for all $X$ and $i$,
$$ CH_0(X,i)_{\Q_l}\oplus CH_0(X,i+1)_{\Q_l}\to H_i(X_\et,\hat \Q_l).$$
Finally, conjecture $P(0)$ can also be recovered from, and
implies a structure theorem for higher Chow groups of
smooth affine schemes: For all smooth and affine schemes $U$
of dimension $d$ over $\F_q$, the groups $CH_0(U,i)$ are
torsion for $i\not= d$, and the canonical map
$CH_0(U,d)_{\Q_l}\to H_d(\bar U_\et,\hat \Q_l)^{Gal(\F_q)}$
is an isomorphism.

Finally, we reproduce an argument of Levine showing that if $F$ is
the absolute Frobenius,
the push-forward $F_*$ acts like $q^n$ on $CH_n(X,i)$, and the pull-back
$F^*$ acts on motivic cohomology $H^i_M(X,\Z(n))$ like $q^n$ for all $n$.
As a Corollary, Conjecture $P(0)$ follows from
finite dimensionality of smooth and projective schemes over finite fields
in the sense of Kimura \cite{kimura}.

\medskip
{\it Acknowledgements:}  This paper was inspired by the work of, and
discussions with, U.Jannsen and S.Saito. We are indebted to the referee,
whose careful reading helped to improve the exposition.

\section{Parshin's conjecture}
We fix a perfect field $k$ of characteristic $p$, and consider
the category of separated schemes of finite type over $k$.
We recall some facts on Bloch's higher Chow groups \cite{bloch},
see \cite{ichhandbook} for a survey. Let $z_n(X,i)$ be the free abelian
group generated by cycles of
dimension $n+i$ on $X\times_k \Delta^i$ which meet all faces
properly, and let $z_n(X,*)$ be the complex of abelian groups obtained
by taking the alternating sum of intersection with face maps as differential.
We define $CH_n(X,i)$ as the $i$th homology of this
complex and motivic Borel-Moore homology to be 
$$H_i^\bm(X,\Z(n))=CH_n(X,i-2n).$$
For a proper map $f:X\to Y$ we have a push-forward
$z_n(X,*)\to z_n(Y,*)$, for a flat, quasi-finite map $f:X\to Y$, we have
a pull-back $z_n(X,*)\to z_n(Y,*)$, and a closed embedding $i:Z\to X$
with open complement $j:U\to X$ induces a localization sequence
$$ \cdots \to H_i^\bm(Z,\Z(n)) \xrightarrow{i_*} H_i^\bm(X,\Z(n))
\xrightarrow{j^*} H_i^\bm(U,\Z(n)) \to\cdots .$$
If $X$ is smooth of pure dimension $d$, then
$H_i^\bm(X,\Z(n))\cong H^{2d-i}(X,\Z(d-n))$, where the right hand side
is Voevodsky's motivic cohomology \cite{voevodsky}.
For a finitely generated field $F$ over $k$, we define
$H_i^\bm(F,\Z(n))=\colim_U H_i^\bm(U,\Z(n))$,
where the colimit runs
through $U$ of finite type over $k$ with field of functions $F$.
For the reader who is more familiar with motivic cohomology,
we mention that Voevodsky's theorem implies that for a field $F$ of
transcendence degree $d$ over $k$, we have
$$H_i^\bm(F,\Z(n))\cong H^{2d-i}(F,\Z(d-n))\cong
\begin{cases}
0&i<d+n\\
K_{d-n}^M(F)&i=d+n.
\end{cases}$$
The latter isomorphism is due to Nesterenko-Suslin and Totaro.
It follows formally from localization that there are spectral sequences
\begin{equation}\label{niveau}
E^1_{s,t}(X)=\bigoplus_{x\in X_{(s)}}H_{s+t}^\bm(k(x),\Z(n))
\Rightarrow H_{s+t}^\bm(X,\Z(n)).
\end{equation}
Here $X_{(s)}$ denotes points of $X$ of dimension $s$.
Since $H_i^\bm(F,\Z(n))=0$ for $i<n+\trdeg F$, the spectral
sequence is concentrated in the area $0\leq s\leq \dim X$ and $t\geq n$.
If we let
$$\tilde H_i^\bm(X,\Z(n))=E^2_{i+n,n}(X)$$ 
be the $i$th homology of the complex
\begin{equation}\label{rational}
0\leftarrow \bigoplus_{x\in X_{(n)}}H_{2n}^\bm(k(x),\Z(n))\leftarrow
\cdots \leftarrow \bigoplus_{x\in X_{(s)}}H_{s+n}^\bm(k(x),\Z(n))
\leftarrow\cdots,
\end{equation}
with $\oplus_{x\in X_{(s)}}H_{s+n}^\bm(k(x),\Z(n))$ in degree $s+n$,
then we obtain a canonical and functorial map
$$\alpha: H_i^\bm(X,\Z(n))\to \tilde H_i^\bm(X,\Z(n)).$$
Note that the groups in \eqref{rational} are Milnor-K-groups.

Parshin's conjecture states that for all smooth and projective
$X$ over $\F_q$, the groups $K_i(X)_\Q$ are torsion for $i>0$.
If Tate's conjecture holds and rational equivalence and homological
equivalence agree up to torsion for all $X$,
then Parshin's conjecture holds by \cite{ichtate}.
Since $K_i(X)_\Q=\oplus_n CH_n(X,i)_\Q$, it follows that Parshin's
conjecture is equivalent to the following conjecture for all $n$.

\medskip

\noindent{\bf Conjecture $P(n)$:}
{\it For all smooth and projective schemes $X$ over the finite field $\F_q$,
the groups $H_i^\bm(X,\Q(n))$ vanish for $i\not=2n$.}

\medskip

We will refer to the following equivalent statements as Conjecture $A(n)$:

\begin{proposition}\label{Ha}
For a fixed finite field $\F_q$ and integer $n$, 
the following statements are equivalent:
\begin{enumerate}
\item For all schemes $X/\F_q$ and all $i$, $\alpha$ induces an
isomorphism $H_i^\bm(X,\Q(n))\cong \tilde H_i^\bm(X,\Q(n))$.
\item For all finitely generated fields $k/\F_q$
with $d:=\trdeg k/\F_q$, and all $i\not=d+n$, we have $H_i^\bm(k,\Q(n))=0$.
\item For all smooth and projective $X$ over $\F_q$ and all $i>\dim X+n$,
we have $H_i^\bm(X,\Q(n))=0$.
\item For all smooth and affine schemes $U$ over $\F_q$ and all
$i>\dim U+n$, we have $H_i^\bm(U,\Q(n))=0$.
\end{enumerate}
\end{proposition}

\proof
\noindent a) $\Rightarrow$ c), d): The complex \eqref{rational}
is concentrated in degrees $[2n,d+n]$.

\noindent c) $\Rightarrow$ b):
This is proved by induction on the transcendence degree.
By de Jong's theorem, we find a smooth and proper model $X$
of a finite extension of $k$. Looking at the niveau spectral
sequence \eqref{niveau}, we see that the induction hypothesis
implies $H_i^\bm(X,\Q(n))=0$ for $i>d+n$
(see \cite{ichtate} for details).

\noindent d) $\Rightarrow$ b): follows by writing $k$ as a colimit
of smooth affine scheme schemes of dimension $d$.

\noindent b) $\Rightarrow$ a): The niveau spectral sequence
collapses to the complex \eqref{rational}.
\proofend

Using the Gersten resolution, the statement in the proposition
implies that on a smooth $X$ of dimension $d$, the motivic complex
$\Q(d-n)$ is concentrated in degree $d-n$, and if
${\mathcal C}_n:={\mathcal H}^{d-n}(\Q(d-n))=CH_n(-,d-n)_\Q$, then
$CH_n(X,i)=H^{d-n-i}(X,{\mathcal C_n})$.

Statement \ref{Ha}d) is part of the following affine analog of $P(n)$:

\medskip

\noindent{\bf Conjecture $L(n)$:}
{\it For all smooth and affine schemes $U$
of dimension $d$ over the finite field $\F_q$, the group $H_i^\bm(U,\Q(n))$
vanishes unless $d\leq i\leq d+n$.}

\medskip

Since $H_i^\bm(U,\Q(n))\cong H^{2d-i}(U,\Q(d-n))$, Conjecture
$L(n)$ can be thought of as an analog of the affine Lefschetz theorem.

\section{Weight homology}
This section is inspired by Jannsen \cite{jannsen}. Throughout this
section we assume resolution of singularities over the field $k$.
Let $\C$ be the category with objects smooth projective varieties over a
field $k$ of characterisic $0$, 
and $\Hom_\C(X,Y)=\oplus_iCH^{\dim Y_i}(X\times Y_i)$,
where $Y_i$ runs through the connected components of $Y$.
Let $\cal H$ be the homotopy category of bounded complexes over $\C$.
Gillet and Soul\'e define in \cite{gilletsoule}, for every separated
scheme of finite type, a
weight complex $W(X)\in \cal H$ satisfying the following properties
\cite[Thm. 2]{gilletsoule} (our notation differs from loc.cit. in variance):

\begin{enumerate}
\item $W(X)$ is represented by a bounded complex
$$ M(X_0)\gets M(X_1)\gets \dots \gets M(X_k)$$
with $\dim X_i\leq \dim X-i$, where $M(X_i)$ placed in degree $i$.
\item $W(-)$ is covariant functorial for proper maps.
\item $W(-)$ is contravariant functorial for open embeddings.
\item If $T\to X$ is a closed embedding with open complement $U$,
then there is a distinguished triangle
$$ W(T)\stackrel{i_*}{\longrightarrow} W(X) \stackrel{j^*}{\longrightarrow}
W(U).$$
\item If $D$ is a divisor with normal crossings in a scheme $X$,
with irreducible components $Y_i$, and if
$Y^{(r)}=\coprod_{\#I=r}\cap_{i\in I}Y_i$, then
$W(X-D)$ is represented by
$$M(X)\gets M(Y^{(1)})\gets \cdots \gets M(Y^{(\dim X)}).$$
\end{enumerate}

The argument in loc.cit. only uses that resolution of singularities
exists over $k$, and we assume from now on that $k$ is such a field.

Given an additive covariant functor $F$ from $\C$ to an abelian category,
we define weight (Borel-Moore) homology $H_i^W(X,F)$ as the $i$th homology
of the complex $F(W(X))$. Weight homology has the functorialities
inherited from b) and c), and satisfies a localization sequence
deduced from d).
If $K$ is a finitely generated field over $k$, then we define
$H_i^W(K,F)$ to be $\colim H_i^W(U,F)$,
where the (filtered) limit runs through integral varieties having
$K$ as their function field.
Similarly, a contravariant functor $G$ from $\C$ to an abelian
category gives rise to weight cohomology (with compact support) $H^i_W(X,G)$.

As a special case, we define the weight homology group
$H_i^W(X,\Z(n))$ as the $i-2n$th homology of the
homological complex of abelian groups $CH_n(W(X))$.

\begin{lemma}\label{popoa}
We have $H_i^W(X,\Z(n))=0$ for $i>\dim X+n$. In particular,
$H_i^W(K,\Z(n))=0$ for every finitely generated field $K/k$ and every
$i>\trdeg_k K+n$.
\end{lemma}

\proof
This follows from the first property of
weight complexes together with $CH_n(T)=0$ for $n>\dim T$.
\proofend

It follows from Lemma \ref{popoa} that the niveau spectral sequence
\begin{equation}\label{uuiioo}
E^1_{s,t}(X)=\oplus_{x\in X_{(s)}}H_{s+t}^W(k(x),\Z(n))\Rightarrow
H_{s+t}^W(X,\Z(n))
\end{equation}
is concentrated on and below the line $t=n$. Let
$$\tilde H_i^W(X,\Z(n))=E^2_{i+n,n}(X)$$ 
be the $i$th homology of the complex \begin{equation}\label{jannsenco}
0\leftarrow \oplus_{x\in X_{(n)}}H_{2n}^W(k(x),\Z(n))\leftarrow
\cdots \leftarrow \oplus_{x\in X_{(s)}}H_{s+n}^W(k(x),\Z(n))\leftarrow\cdots,
\end{equation}
where $\oplus_{x\in X_{(s)}}H_{s+n}^W(k(x),\Z(n))$ is placed in degree $s+n$.
Then we obtain a canonical and natural map  
$$\gamma:\tilde H_i^W(X,\Z(n))\to H_i^W(X,\Z(n)).$$

Consider the canonical map of covariant functors $\pi':z_n(-,*)\to CH_n(-)$
on the category of smooth projective schemes over $k$,
sending the cycle complex
to its highest cohomology. Then by \cite[Thm.5.13, Rem.5.15]{jannsen},
the set of associated homology functors extends to a homology theory
on the category of allvarieties over $k$. 
The argument of \cite[Prop. 5.16]{jannsen}
show that the extension of the associated homology functors for $z_n(-,*)$
are higher Chow groups $CH_n(-,i)$. The extension $CH_n(-)$ are by definition
the functors $H_i^W(-,\Z(n))$. We obtain a functorial map
$$
\pi:H_i^\bm(X,\Z(n))\to H_i^W(X,\Z(n)).
$$

\begin{lemma}\label{degreezero}
For $i=2n$, and all schemes $X$, the map $\pi$ is an isomorphism
$H_{2n}^\bm(X,\Z(n))\cong H_{2n}^W(X,\Z(n))$.
In particular, $H_{2d}^W(K,\Z(d))\cong \Z$ for all fields $K$ of
transcendence degree $d$ over $k$.
\end{lemma}

\proof
The statement is clear for $X$ smooth and projective.
We proceed by induction on the dimension of $X$. Using the
localization sequence for both theories, we can assume that
$X$ is proper. Let $f:X'\to X$ be a resolution of singularities of $X$,
$Z$ be the closed subscheme (of lower dimension) where $f$ is not an
isomorphism, and $Z'=Z\times_XX'$.
Then we conclude by comparing localization sequences
$$\begin{CD}
H_{2n}^\bm(Z',\Z(n))@>>> H_{2n}^\bm(Z,\Z(n))\oplus H_{2n}^\bm(X',\Z(n))
@>>>H_{2n}^\bm(X,\Z(n)) @>>>0\\
@| @| @VVV \\
H_{2n}^W(Z',\Z(n))@>>> H_{2n}^W(Z,\Z(n))\oplus H_{2n}^W(X',\Z(n))@>>>
H_{2n}^W(X,\Z(n)) @>>>0.
\end{CD}$$
\proofend

The map $\pi$ for fields induces a map
$\beta:\tilde H_i^\bm(X,\Q(n))\to \tilde H_i^W(X,\Q(n))$,
which fits into the (non-commutative) diagram
\begin{equation*}
\begin{CD}
H_i^\bm(X,\Z(n))@>\pi>> H_i^W(X,\Z(n))\\
@V\alpha VV @AA\gamma A\\
\tilde H_i^\bm(X,\Z(n))@>\beta>> \tilde H_i^W(X,\Z(n)).
\end{CD}
\end{equation*}

We now return to the situation $k=\F_q$, and compare weight homology to
higher Chow groups using their niveau spectral sequences. We saw
that the niveau spectral sequence for higher Chow groups is concentrated
above the line $t=n$, and that the niveau spectral sequence for weight
homology is concentrated below the line $t=n$. Our aim is to show
that Parshin's conjecture is
equivalent to both being rationally concentrated on this line, and that
the resulting complexes are isomorphic.

The following statements will be referred to as Conjecture $B(n)$:

\begin{proposition}\label{Hc}
For a fixed integer $n$, the following statements are equivalent:
\begin{enumerate}
\item The map $\beta$ induces an isomorphism
$\tilde H_i(X,\Q(n))\cong \tilde H_i^W(X,\Q(n))$
for all schemes $X$ and all $i$.
\item The map $\pi$ induces an isomorphism
$H_{d+n}^\bm(k,\Q(n))\cong H_{d+n}^W(k,\Q(n))$
for all finitely generated fields $k/\F_q$, where $d=\trdeg k/\F_q$.
\item For every smooth and projective $X$ over $\F_q$
we have $\tilde H_{d+n}^\bm(X,\Q(n))=\tilde H_{d+n-1}^\bm(X,\Q(n))=0$
if $d=\dim X>n+1$, and $\tilde H_{2n+1}^\bm(X,\Q(n))=0$ if $\dim X=n+1$.
\end{enumerate}
\end{proposition}

Note that assuming $A(n)$, c) is equivalent to
$H_{d+n}^\bm(X,\Q(n))=H_{d+n-1}^\bm(X,\Q(n))=0$ and $H_{2n+1}^\bm(X,\Q(n))=0$
for all smooth and projective
$X$ of dimension $d>n+1$ and $d=n+1$, respectively,
hence are part of Conjecture $P(n)$.

\medskip

\proof
\noindent b) $\Rightarrow$ a) is trivial, and a) $\Rightarrow$ b)
follows by a colimit argument because
$$\colim_{U\subseteq X} \tilde H_i^\bm(X,\Q(n))\cong
\begin{cases}
H_{d+n}^\bm(k(X),\Q(n)) &i=d+n;\\
0 &\text{otherwise.}
\end{cases}$$

\noindent a) $\Rightarrow$ c) follows because for $X$ smooth and
projective of dimension $d$, the cohomology of the complex
\eqref{rational} tensored with $\Q$ equals
$\tilde H_i^\bm(X,\Q(n))=\tilde H^W_i(X,\Q(n))$
for $i=d+n$ and $i=d+n-1$ (or $i=2n+1$ in case $d=n+1$).
An inspection of the niveau spectral sequence \eqref{uuiioo}
shows that this is a subgroup of $H^W_i(X,\Q(n))=0$.

\noindent c) $\Rightarrow$ b): For $n>d$, both sides vanish,
whereas for $d=n$, both sides are canonically isomorphic to $\Q$.
For $n<d$, we proceed by induction on $d$.
Choose a smooth and projective model $X$ for $k$
and compare the exact sequences \eqref{rational} and \eqref{jannsenco}
$$\begin{CD}
A @<<<\oplus_{x\in X_{(d-1)}}H_{d+n-1}^\bm(k(x),\Q(n))
@<<< H_{d+n}^\bm(k,\Q(n))\\
@| @|@VVV \\
B @<<<\oplus_{x\in X_{(d-1)}}H_{d+n-1}^W(k(x),\Q(n))
@<<< H_{d+n}^W(k,\Q(n))
\end{CD}$$
The terms on the left are $A=H_{2n}^\bm(X,\Q(n)) \cong B=H_{2n}^W(X,\Q(n))$
if $d=n+1$, and $A=\oplus_{x\in X_{(d-2)}}H_{d+n-2}^\bm(k(x),\Q(n))\cong
B=\oplus_{x\in X_{(d-2)}}H_{d+n-2}^W(k(x),\Q(n))$, if $d>n+1$.
The upper sequence is exact by hypothesis, and an inspection of \eqref{uuiioo}
shows that the lower sequence is exact because $H^W_i(X,\Q(n))=0$ for $i>2n$.
\proofend

We refer to the following statements as Conjecture C(n):

\begin{proposition}\label{Hb}
For a fixed integer $n$, the following statements are equivalent:
\begin{enumerate}
\item For all schemes $X$ over $\F_q$ and all $i$, the map
$\gamma$ induces an isomorphism $\tilde H_i^W(X,\Q(n)\cong H_i^W(X,\Q(n))$.
\item For all finitely generated fields $k/\F_q$ and $i\not=\trdeg k/\F_q+n$,
we have $H_i^W(k,\Q(n))=0$.
\item For all smooth and projective $X$, the map $\gamma$ induces an isomorphism
$$\tilde H_i^W(X,\Q(n))=
\begin{cases}
0&i>0;\\
CH_n(X)_\Q&i=2n.
\end{cases}$$
\end{enumerate}

\end{proposition}

\proof
\noindent b) $\Rightarrow$ a) is trivial and  a) $\Rightarrow$ b)
follows by a colimit argument.

\noindent a) $\Rightarrow$ c) is trivial for $i>2n$,
and Lemma \ref{degreezero} for $i=2n$.

\noindent c) $\Rightarrow$ b):
This is proved like Proposition \ref{Ha},
by induction on the transcendence degree of $k$.
Let $X$ be a smooth and projective model of $k$.
The induction hypothesis implies that the niveau spectral sequence
\eqref{uuiioo} collapses to the horizontal line $t=n$ and the
vertical line $s=d$.
Since it converges to $H_i^W(X,\Q(n))$, which is zero for $i>2n$,
we obtain isomorphisms
$H_{d+n-i+1}^W(k(X),\Q(n)) \xrightarrow{d_i}\tilde H_{d+n-i}^W(X,\Q(n))$
for $1<i<d-n$, and an exact sequence
$$H_{2n+1}^W(k(X),\Q(n))\stackrel{d_{d-n}}{\hookrightarrow}
\tilde H_{2n}^W(X,\Q(n))\to H_{2n}^W(X,\Q(n))\twoheadrightarrow H_{2n}^W(k(X),\Q(n)).$$
The claim follows because
$\tilde H_i^W(X,\Q(n))= H_i^W(X,\Q(n))=0$ for $i>2n$, and because
the maps $H_{2n}^\bm(X,\Q(n))\xrightarrow{\pi} H_{2n}^W(X,\Q(n))
\xleftarrow{\gamma}\tilde H_{2n}^W(X,\Q(n))$ are isomorphisms by
Lemma \ref{degreezero} and hypothesis.
\proofend

\begin{theorem}\label{Hall}
The conjunction of Conjectures $A(n)$, $B(n)$ and $C(n)$ is equivalent
to $P(n)$.
\end{theorem}

\proof
Given Conjectures $A(n)$, $B(n)$, and $C(n)$, we get
$$ H_i^\bm(X,\Q(n))=\tilde H_i^\bm(X,\Q(n))\cong \tilde H_i^W(X,\Q(n))
\cong H_i^W(X,\Q(n))$$
for all $X$, and the latter vanishes
for smooth and projective $X$ and $i>0$, hence Conjecture $P(n)$ follows.
Conversely, Conjecture $P(n)$ implies Prop. \ref{Ha}c), then \ref{Hc}c),
and finally \ref{Hb}c) by using \ref{Ha}a) and \ref{Hc}a).
\proofend

\rem
Propositions \ref{Ha}, \ref{Hb}, \ref{Hc} as well as Theorem
\ref{Hall} remain true if we restricts ourselves to schemes of
dimension at most $N$, and fields of transcendence degree at most $N$,
for a fixed integer $N$.

\rem
Gillet announced that one can obtain a rational version of weight
complexes by using de Jong's theorem
on alterations instead of resolution of singularities. The same
argument should then give the generalization \cite[Thm.5.13]{jannsen}.
In this case, all arguments of this section hold true rationally,
except the proof of c) $\Rightarrow$ b) in Propositions \ref{Hc}
and \ref{Hb}, and the proof of $P(n)\Rightarrow A(n), B(n), C(n)$
in Theorem \ref{Hall},
which require that every finitely generated field over $\F_q$
has a smooth and projective model.

\section{The case $n=0$}
Since $H_i^\bm(X,\Q(0))=CH_0(X,i)_\Q$, we use higher Chow groups in
this section.

\begin{proposition}\label{surf}
We have $CH_0(X,i)_\Q\cong \tilde H_i^\bm(X,\Q(0))$ for $i\leq 2$,
and the map $CH_0(X,3)_\Q\to \tilde H_3^\bm(X,\Q(0))$ is surjective
for all $X$. In particular, $A(0)$ holds in dimensions at most $2$.
\end{proposition}

\proof
Since $H^i(k(x),\Q(1))=0$ for $i\not=1$ and $H^i(k(x),\Q(0))=0$ for $i\not=0$,
this follows from an inspection of the niveau spectral sequence.
\proofend

Jannsen \cite{jannsen} defines a variant of weight homology
with coefficients $A$, $H^W_i(X,A)$,  as the
homology of the complex $\Hom(CH^0(W(X)),A)$. Note that
$H^W_i(X,\Q)=H^W_i(X,\Q(0))$ because $\Hom(CH^0(X),\Q)\cong CH_0(X)_\Q$
for smooth and projective. $X$, in a functorial way.
Indeed, a map of connected, smooth and projective varieties induces the
identity pull-back on $CH^0$ and the identity push-forward on $CH_0$.

\begin{theorem}\label{jweight0} (Jannsen)
Under resolution of singularities,
$H_a^W(k,A)=0$ for $a\not=\trdeg k/\F_q$, hence $H_i^W(X,A)$ is the
homology of the complex
\begin{equation}\label{jannsen}
0\leftarrow \oplus_{x\in X_{(0)}}H_0^W(k(x),A)\leftarrow
\cdots \leftarrow \oplus_{x\in X_{(s)}}H_s^W(k(x),A)\leftarrow\cdots
\end{equation}
for all schemes $X$. In particular, Conjecture C(0) holds.
\end{theorem}

This is proved in \cite[Prop. 5.4, Thm. 5.10]{jannsen}.
The proof only works for $n=0$, because it
uses the bijectivity of $CH^0(Y)\to CH^0(X)$ for
a map of connected smooth and projective schemes $X\to Y$.
The second statement follows using the niveau spectral sequence,
which exists because $H^W_i(X,A)$ satisfies the localization
property by property (4) of weight complexes.

\medskip

Let $\Z^c(0)$ be the complex of etale sheaves $z_0(-,*)$.
For any prime $l$, consider $l$-adic cohomology
$$H_i(X_\et,\hat \Q_l):= \Q\otimes_\Z \lim H_i(X_\et,\Z^c/l^r(0)).$$
In \cite{ichdual}, we showed that for every positive integer $m$,
and every scheme $f:X\to k$ over a perfect field, there is a quasi-isomorphism
$\B/m(0)\cong Rf^!\Z/m$. In particular, the above definition
agrees with the usual definition of $l$-adic homology if $l\not=p=\chr \F_q$.
If $\bar X=X\times_{\F_q}\bar \F_q$ and $\hat G=\Gal(\bar \F_q/\F_q)$,
then there is a short exact sequence
$$ 0 \to H_{i+1}(\bar X_\et,\hat \Q_l)_{\hat G} \to H_i(X_\et,\hat \Q_l)
\to H_i(\bar X_\et,\hat \Q_l)^{\hat G}\to 0$$
and for $U$ affine and smooth, $H_i(U_\et,\hat \Q_l)$ vanishes
for $i\not=d, d-1$ and $l\not=p$ by the affine Lefschetz theorem and
a weight argument \cite[Thm.3a)]{kahn}.
The map from Zariski-hypercohomology of $\B/m(0)$ to etale-hypercohomology
of $\B/m(0)$ induces a functorial map
$$ CH_0(X,i)/m\to CH_0(X,i,\Z/m)\to H_i(X_\et,\Z^c/m(0)),$$
hence in the limit a map
$$ \omega:CH_0(X,i)_{\Q_l} \to H_i(X_\et,\hat \Q_l).$$
Similarly, the map
$ \bar\omega:CH_0(\bar X,i)_{\Q_l} \to H_i(\bar X_\et,\hat \Q_l)$
induces a map
$$\tau :CH_0(X,i+1)_{\Q_l}\xleftarrow{\sim} (CH_0(\bar X,i+1)_{\Q_l})_{\hat G}
\to H_{i+1}(\bar X_\et,\hat \Q_l)_{\hat G} \to H_i(X_\et,\hat \Q_l)$$
(the left map is an isomorphism by a trace argument).
The sum
\begin{equation}\label{gamm}
\varphi_X^i:  CH_0(X,i)_{\Q_l}\oplus CH_0(X,i+1)_{\Q_l}\to H_i(X_\et,\hat \Q_l)
\end{equation}
is compatible with localization sequences, because all maps involved
in the definition are.
The following proposition shows that Parshin's conjecture can
be recovered from and implies a structure theorem for
higher Chow groups of smooth affine schemes;
compare to Jannsen \cite[Conj. 12.4b)]{jannsen}.

\begin{proposition}\label{parshineq}
The following statements are equivalent:
\begin{enumerate}
\item Conjecture $P(0)$.
\item For all schemes $X$, and all $i$, the map $\varphi_X^i$
is an isomorphism.
\item For all smooth and affine schemes $U$ of dimension $d$,
the groups $CH_0(U,i)$ are torsion for $i\not= d$, and the composition
$\omega: CH_0(U,d)_{\Q_l}\to H_d(U_\et,\hat \Q_l)\to
H_d(\bar U_\et,\hat \Q_l)^{\hat G}$ is an isomorphism.
\end{enumerate}
\end{proposition}

\proof
\noindent
a) $\Rightarrow$ b): First consider the case that $X$ is smooth
and proper. Then Conjecture $P(0)$ is equivalent to vanishing of the
left hand sides of \eqref{gamm} for $i\not=0,-1$, whereas the right
hand side of \eqref{gamm} vanishes by the Weil-conjectures.
On the other hand, $\varphi_X^0$ induces an isomorphism
$CH_0(X)\otimes\Q_l\cong H^{2d}(X_\et,\hat \Q_l(d))$
and $\varphi_X^{-1}$ induces an isomorphism
$CH_0(X)\otimes \Q_l \cong (CH_0(\bar X)\otimes \Q_l)_{\hat G}
\cong H^{2d}(\bar X,\hat \Q_l(d))_{\hat G}$. Indeed, both sides
are isomorphic to $\Q_l$ if $X$ is connected. Using localization and
alterations, the statement for smooth and proper $X$ implies the
statement for all $X$.

\noindent b) $\Rightarrow$ a): The right-hand side of \eqref{gamm} is zero
for $i\not=0,-1$ by weight reasons for smooth and projective $X$,
hence so is the left side.

\noindent
b) $\Rightarrow$ c): This follows because $H_i(U_\et,\hat \Q_l)=0$
unless $i=d,d-1$ for smooth and affine $U$.

\noindent
c) $\Rightarrow$ b): We first assume that $X$ is smooth and affine.
By hypothesis and the affine Lefschetz
theorem, both sides of \eqref{gamm} vanish for $i\not=d,d-1$, and are
isomorphic for $i=d$. For $i=d-1$, the vertical maps in
the following diagram are isomorphisms by semi-simplicity,
$$\begin{CD}
CH_0(X,d)_{\Q_l}\cong CH_0(\bar X,d)_{\Q_l}^{\hat G}@>>>
H_d(\bar X,\hat \Q_l)^{\hat G} \\
@| @| \\
(CH_0(\bar X,d)_{\Q_l})_{\hat G}@>>>
H_d(\bar X,\hat \Q_l)_{\hat G}@>\sim >> H_{d-1}(X,\hat \Q_l).
\end{CD}$$
Hence the lower map $\varphi_X^{d-1}$ is an isomorphism because
the upper map is.
Using localization, the statement for smooth and affine $X$ implies the
statement for all $X$.
\proofend

\begin{proposition}
Under resolution of singularities, the following are equivalent:
\begin{enumerate}
\item Conjecture $P(0)$.
\item For every smooth affine $U$
of dimension $d$ over $\F_q$, we have $CH_0(U,i)_\Q\cong H_i^W(U,\Q)$
for all $i$, and these group vanish for $i\not=d$.
\item For every smooth affine $U$ of dimension $d$ over $\F_q$,
the groups $CH_0(U,i)_\Q$ vanish for $i>d$, and
$CH_0(U,d)_\Q\cong H_d^W(U,\Q)$.
\end{enumerate}
\end{proposition}

\proof
a) $\Rightarrow$ b): It follows from the previous Proposition that
$CH_0(U,i)_\Q=0$ for $i\not=d$. On the other hand, $P(0)$
for all $X$ implies that $CH_0(X,i)\cong H_i^W(X,\Q)$ for all $i$
and $X$.

c) $\Rightarrow$ a): The statement implies Conjecture $A(0)$, and
then Conjecture $B(0)$, version b), for all $X$. By Theorems \ref{Hall}
and \ref{jweight0}, $P(0)$ follows.
\proofend

\begin{proposition}
Conjecture $P(0)$ for all smooth and projective $X$ implies
the following statements:

a) (Affine Gersten) For every smooth affine $U$ of dimension $d$,
the following sequence is exact:
$$ CH_0(U,d)_\Q \hookrightarrow \oplus_{x\in U^{(0)}}H^d(k(x),\Q(d))
\to \oplus_{x\in U^{(1)}}H^{d-1}(k(x),\Q(d-1))\to \cdots .$$

b) Let $X=X_d\supseteq X_{d-1}\supseteq \cdots X_1\supseteq X_0$
be a filtration such that $U_i=X_i-X_{i-1}$ is smooth and affine
of dimension $i$.
Then $CH_0(X,i)_\Q$ is isomorphic to the $i$th homology of the complex
$$ 0\to CH_0(U_d,d)_\Q \to CH_0(U_{d-1},d-1)_\Q\to \cdots
\to CH_0(U_0,0)_\Q\to 0.$$
The maps
$CH_0(U_i,i)_\Q\to CH_0(X_{i-1},i-1)_\Q\to CH_0(U_{i-1},i-1)_\Q$
arise from the localization sequence.
\end{proposition}

\proof
a) follows because the spectral sequence
$\eqref{niveau}$ collapses, and b) by a diagram chase.
\proofend

\rem If we fix a smooth scheme $X$ of dimension $d$, and use cohomological
notation, then by Proposition \ref{Ha} and the Gersten resolution, we get
that the rational motivic complex $\Q(d)$ is conjecturally concentrated in
degree $d$, say $\C_d={\mathcal H}^d(\Q(d))=CH_0(-,d)_\Q=H_d^W(-,\Q)$.
Then Conjecture $L(0)$ says that $H^i(U,\C_d)=0$ for $U\subseteq X$ affine
and $i>0$. This is analog to the mod $p$ situation, where
the motivic complex agrees with the logarithmic de Rham-Witt sheaf
$\Z/p(n)\cong \nu^d[-d]$, and $H^i(U_\et,\nu^d)=0$ for $U\subseteq X$ affine
and $i>0$. The latter can be proved by writing $\nu^d$ as the kernel
of a map of coherent sheaves and using the vanishing of cohomology
of coherent sheaves on affine schemes. This suggest that one might
try to do the same for $\C_d$.

\subsection{Frobenius action}
Let $F:X\to X$ be the Frobenius morphism induced by the $q$th
power map on the structure sheaf.

\begin{theorem}
The push-forward $F_*$ acts like $q^n$ on $CH_n(X,i)$, and the pull-back
$F^*$ acts on $H^i(X,\Z(n))$ as $q^n$ for all $n$.
\end{theorem}

The Theorem is well-known, but we could not find a proof in the literature.
The proof of Soul\'e \cite[Prop.2]{soule} for Chow groups does not carry
over to higher Chow groups, because the Frobenius does not
act on the simplices $\Delta^n$, hence a cycle $Z\subseteq \Delta^n\times X$
is not send to a multiple of itself by the Frobenius.
We give an argument due to M. Levine.

\proof
Let $DM^-$ be Voevodsky's derived category of bounded above
complexes of Nisnevich sheaves with transfers with homotopy
invariant cohomology sheaves. Then we have the isomorphisms
\begin{align*}
CH_n(X,i)&\cong \Hom_{DM^-}(\Z(n)[2n+i],M_c(X)),\\
H^i(X,\Z(n))&\cong \Hom_{DM^-}(M(X),\Z(n)[i]).
\end{align*}
The action of
the Frobenius is given by composition with $F:M_c(X)\to M_c(X)$
and $F:M(X)\to M(X)$, respectively.

The Frobenius acts on the category $DM^-$, i.e. for
every $\alpha\in \Hom_{DM^-}(X,Y)$ we have $F_Y\circ\alpha=\alpha\circ F_X$.
This follows by considering composition of correspondences.
Hence it suffices to calculate the action of the Frobenius on $\Z(n)$,
i.e. show that $F=q^n\in \Hom_{DM^-}(\Z(n),\Z(n))\cong \Z$. But
$\Hom_{DM^-}(\Z(n),\Z(n))$
is a direct factor of $\Hom_{DM^-}(\Z(n),{\mathbb P}^n[-2n])=
CH_n({\mathbb P}^n)$. The latter is the free abelian group generated
by the generic point, and the Frobenius acts by $q^n$ on it.
\proofend

\rem 1) It would be interesting to write down an explicit chain homotopy
between $F_*$ and $q^n$ on $z_n(X,*)$.

2) The proposition implies that the groups $CH_n(\F_q,i)$
are killed by $q^n-1$, and that $CH_n(X,i)$ is $q$-divisible for $n<0$.

Granted the Theorem, the standard argument gives the following
Corollary, see also Jannsen \cite[Thm. 12.5.7]{jannsenmurre}.

\begin{corollary}
Assume that for all be smooth and projective $X$ of dimension $d$, the
kernel of the map $CH^d(X\times X)\to \End_{hom}(M(X))$
is nilpotent. Then Conjecture $P(0)$ holds.
\end{corollary}

The hypothesis of the Corollary is satisfied if $X$ is finite
dimensional in the sense of Kimura \cite{kimura}.

\proof
Using the existence of a zero-cycle $c$ of degree $1$, we see
that the projector $\pi_{2d}=[X\times c]$ is defined.
Let $\tilde X$ be the motive $\ker \pi_{2d}=X/{\mathcal L}^d$,
where ${\mathcal L}$ is the Lefschetz motive.
Consider the action of the geometric Frobenius
$F\in \End_{rat}(\tilde X)\subseteq CH^d(X\times X)$.
Its image in the category of motives for homological equivalence is
algebraic, and its minimal polynomial $P_{\tilde X}(T)$ has roots of
absolute value $q^\frac{j}{2}$ for $0\leq j<2d$. By hypothesis,
$P_{\tilde X}(F)^a=0$ in $\End_{rat}(\tilde X)$ for some integer $a$, but
by the Theorem,
$F^*$ acts like $q^d$ on $CH^d(X,i)$. Hence
$0\not=P_{\tilde X}(q^d)^a=P_{\tilde X}(F^*)^a=0$.
\proofend


\begin{thebibliography}{99}
\bibitem{bloch} {\sc S.Bloch}, Algebraic cycles and higher $K$-theory.
Adv. in Math. 61 (1986), no. 3, 267--304.

\bibitem{ichtate} {\sc T.Geisser}, Tate's conjecture, algebraic cycles and
rational $K$-theory in characteristic $p$.
$K$-Theory 13 (1998), no. 2, 109--122.

\bibitem{ichhandbook} {\sc T.Geisser}, Motivic cohomology, K-theory,
and topological cyclic homology.
Handbook of $K$-theory. Vol. 1, 2, 193--234, Springer, Berlin, 2005.

\bibitem{ichdual} {\sc T.Geisser}, Duality via cycle complexes,
Preprint 2006.

\bibitem{ichkato} {\sc T.Geisser}, Arithmetic homology and an integral version
of Kato's conjecture. Preprint 2007.

\bibitem{marcI} {\sc T.Geisser, M.Levine}, The $p$-part of
$K$-theory of fields in characteristic $p$. Inv. Math. {\bf 139}
(2000), 459--494.

\bibitem{gillet}{\sc H.Gillet}, Homological descent for the $K$-theory
of coherent sheaves. Algebraic $K$-theory, number theory, geometry
and analysis (Bielefeld, 1982), 80--103, Lecture Notes in Math., 1046,
Springer, Berlin, 1984.

\bibitem{gilletsoule} {\sc H.Gillet, C.Soul\'e}, Descent, motives and
$K$-theory. J. Reine Angew. Math. 478 (1996), 127--176.

\bibitem{jannsenln} {\sc U.Jannsen}, Mixed motives and algebraic $K$-theory.
Lecture Notes in Mathematics, 1400. Springer-Verlag, Berlin, 1990.

\bibitem{jannsenmurre} {\sc U.Jannsen}, Some remarks on finite
dimensional motives and Murre's conjecture, (2006).

\bibitem{jannsen} {\sc U.Jannsen}, Hasse principles for higher-dimensional
fields, Preprint Universit\"at Regensburg 18/2004.

\bibitem{kahn} {\sc B.Kahn}, Some finiteness results for
etale cohomology. J. Number Theory 99 (2003), no. 1, 57--73.

\bibitem{kimura} {\sc S.Kimura}, Chow groups are finite dimensional,
in some sense. Math. Ann. 331 (2005), no. 1, 173--201.

\bibitem{soule} {\sc C.Soule}, Groupes de Chow et $K$-theorie de varietes
sur un corps fini. Math. Ann. 268 (1984), no. 3, 317--345.

\bibitem{voevodsky} {\sc V.Voevodsky}, Motivic cohomology groups are
isomorphic to higher Chow groups in any characteristic.
Int. Math. Res. Not. 2002, no. 7, 351--355.

\end{thebibliography}
\end{document}